\documentstyle[amssymb,12pt]{amsart}
\newtheorem{prop}{Proposition}[section]

\newtheorem{lemma}{Lemma}[section]
\newtheorem{thm}{Theorem}[section]
\newtheorem{corollary}{Corollary}[section]

\theoremstyle{remark}
\newtheorem{remark}{Remark}

\begin{document}
\newcommand{\nc}{\newcommand}
\nc{\on}{\operatorname}
\nc{\pa}{\partial}
\nc{\cA}{{\cal A}}\nc{\cB}{{\cal B}}\nc{\cC}{{\cal C}}
\nc{\cE}{{\cal E}}\nc{\cG}{{\cal G}}\nc{\cH}{{\cal H}}
\nc{\cX}{{\cal X}}\nc{\cR}{{\cal R}}\nc{\cL}{{\cal L}}
\nc{\cZ}{{\cal Z}}\nc{\cT}{{\cal T}}\nc{\cK}{{\cal K}}
\nc{\sh}{\on{sh}}\nc{\Id}{\on{Id}}
\nc{\ad}{\on{ad}}\nc{\Der}{\on{Der}}\nc{\End}{\on{End}}\nc{\res}{\on{res}}
\nc{\Imm}{\on{Im}}\nc{\limm}{\on{lim}}\nc{\Ad}{\on{Ad}}
\nc{\Hol}{\on{Hol}}\nc{\Det}{\on{Det}}\nc{\tr}{\on{tr}}\nc{\Yg}{\on{Yg}}
\nc{\ev}{\on{ev}}\nc{\Aut}{\on{Aut}}
\nc{\de}{\delta}\nc{\si}{\sigma}\nc{\ve}{\varepsilon}
\nc{\al}{\alpha}
\nc{\CC}{{\Bbb C}}\nc{\ZZ}{{\Bbb Z}}\nc{\NN}{{\Bbb N}}
\nc{\zz}{{\bold z}}\nc{\ww}{{\underline w}}
\nc{\eeta}{{\underline{\eta}}}\nc{\xxi}{{\underline
{\xi}}}\nc{\uu}{{\underline u}}\nc{\xx}{{\underline x}}
\nc{\yy}{{\underline y}}
\nc{\AAA}{{\Bbb A}}\nc{\VV}{{\Bbb V}}
\nc{\cO}{{\cal O}} \nc{\cF}{{\cal F}}
\nc{\la}{{\lambda}}
\nc{\lla}{{\bold \la}}
\nc{\G}{{\frak g}}\nc{\A}{{\frak a}}
\nc{\HH}{{\frak h}}
\nc{\N}{{\frak n}}\nc{\B}{{\frak b}}
\nc{\La}{\Lambda}
\nc{\g}{\gamma}\nc{\eps}{\epsilon}\nc{\wt}{\widetilde}
\nc{\wh}{\widehat}
\nc{\bn}{\begin{equation}}\nc{\en}{\end{equation}}
\nc{\SL}{{\frak{sl}}}\nc{\GL}{{\frak{gl}}}

%
%
%

\newcommand{\ldar}[1]{\begin{picture}(10,50)(-5,-25)
\put(0,25){\vector(0,-1){50}}
\put(5,0){\mbox{$#1$}} 
\end{picture}}

\newcommand{\lrar}[1]{\begin{picture}(50,10)(-25,-5)
\put(-25,0){\vector(1,0){50}}
\put(0,5){\makebox(0,0)[b]{\mbox{$#1$}}}
\end{picture}}

\newcommand{\luar}[1]{\begin{picture}(10,50)(-5,-25)
\put(0,-25){\vector(0,1){50}}
\put(5,0){\mbox{$#1$}}
\end{picture}}

\title[Coinvariants for Yangian doubles and quantum KZ equations]
{Coinvariants for Yangian doubles and quantum Knizhnik-Zamolodchikov 
equations}

\author{B. Enriquez}\address{Centre de Math\'ematiques, URA 169 du CNRS,
Ecole Polytechnique,
91128 Palaiseau, France}

\author{G. Felder}
\address{D-Math, ETH-Zentrum, HG G46, CH-8092 Z\"urich, Suisse}

\date{July 1997}

\begin{abstract} 
We present a quantum version of the construction of the KZ system of equations
as a flat connection on the spaces of coinvariants of representations of
tensor products of Kac-Moody algebras. We consider here representations
of a tensor product of Yangian doubles and compute the coinvariants of 
a deformation of the subalgebra generated by the regular
functions of a rational curve with marked points. 
We observe that Drinfeld's quantum Casimir element can be viewed as a
deformation of the zero-mode of the Sugawara tensor in the
Yangian double. These ingredients serve to define a compatible system of
difference equations, which we identify with the quantum KZ equations
introduced by I. Frenkel and N. Reshetikhin. 
\end{abstract}

\maketitle

\subsection*{Introduction}

The Knizhnik-Zamolodchikov (KZ) system is a set of  
differential equations satisfied by
correlation functions in Wess-Zumino-Witten theories (\cite{KZ}). These
equations can be interpreted as the equations satisfied by matrix
elements of intertwining operators associated with representations of
affine Kac-Moody algebras (\cite{TK}). They define a local system on the
configuration space of $n$ distinct marked points on the rational curve
$\CC P^{1}$. This local system also has another
interpretation: to a complex curve with marked points, a system of
weights of a semisimple Lie algebra and a positive integer is
associated the vector space of conformal blocks. It is the space of
coinvariants of a representation of a product of Kac-Moody algebras
attached to the points of the curve, with respect to the subalgebra
formed by the rational functions on the curve, regular outside the points. 
The conformal blocks form
a vector bundle on the moduli space of curves with marked points. 
There exists a natural connection on this vector bundle, the KZB (for
KZ-Bernard) connection. It is provided
by the action of the Sugawara field. This field is a generating
functional for elements of the
enveloping algebra of the Kac-Moody algebra. One shows that certain
combinations of these elements conjugate the regular
algebras associated to nearby elements of the moduli space (see
\cite{Ber,TUY}). 

In \cite {FR}, I. Frenkel and N. Reshetikhin introduced $q$-deformed
analogues of the KZ equations. This qKZ system is a difference system
obeyed by matrix elements of intertwining operators of quantum affine
algebras. Later, elliptic analogues of this system were defined and
studied (\cite{F,FTV}). If one wished to understand $q$-deformed
versions of the KZB connection in higher genus, it would be important to
understand how these equations could be derived from the coinvariants
viewpoint. To make such a derivation explicit in the rational case
is the main goal of this paper. 

Let us now present our work. 
We consider a system of points $\zz =
(z_{i})_{i=1,\ldots,n}$ on $\CC P^{1}$. We call $\cO_{i}$ and $\cK_{i}$
the completed ring and field of $\CC P^{1}$ at $z_{i}$. We also call
$\cO$ and $\cK$ their direct sums, and $R_{\zz}$ the ring of regular
functions on $\CC P^{1}$, regular  outside $\zz$ and vanishing at
infinity. We set $\bar\G =\SL_{2}$, and we denote by $\G_{\cK}$ the
double extension of $\bar\G \otimes \cK$ by central and derivation
elements, and by $\G_{\cO}$ and
$\G_{\cK}$ extensions $\bar\G\otimes \cO$ and of $\bar\G \otimes
R_{\zz}$. In this situation, the coinvariants construction described
above is
based on the inclusions of the enveloping algebras $U\G_{\cO}$ and
$U\G_{\zz}$ in $U\G_{\cK}$. 

To construct a deformation of these inclusions, we note that the
decomposition of $\G_{\cK}$ as a direct sum $\G_{\cO} \oplus \G_{\zz}$
is that of a Manin triple, associated with the rational form $dz$ on
$\CC P^{1}$. We apply techniques of quantum currents and twists 
\cite{ER2,EF} to the quantization of this triple (sect. \ref{quant}). 
We then give an
presentation of the resulting algebras $U_{\hbar}\G_{\cO},
U_{\hbar}\G_{\zz}$ and $U_{\hbar}\G_{\cK}$ in terms of $L$-operators
(sect. \ref{rll}). 
Let us mention here
that P. Etingof and D. Kazhdan obtained in \cite{EK} quantizations of
these triples, for arbitrary $\bar\G$; their construction 
of $U_{\hbar}\G_{\cK,\zz}$ from the double
Yangian algebra is a special case of their construction of 
``factored algebras'', that applies to any quantum group based on an
$1$-dimensional algebraic group. They also
obtain $RTT$-type relations similar to our $RLL$ relations (\cite{EK},
Prop. 3.25). 
 
Our next step is to construct an isomorphism between $U_{\hbar}\G_{\cK}$
and a tensor product of $n$ copies of the double Yangian algebra
$DY(\SL_{2})$, with their central elements identified. 
This is done in Prop. \ref{isom} using the $L$-operators
description of both algebras. The formulas closely resemble formulas for
coproducts. They may have the following interpretation: the isomorphism
of the quantum groups appearing in \cite{ER1} with a tensor product of
``local'' algebras could be obtained if we extended the construction of
that paper to a larger extension, with a $n$-dimensional center and $n$
derivations. The resulting algebra would then have specialization
morphisms to local algebras, and the desired isomorphism would result
from composing the coproduct $\Delta^{(n)}$ with the tensor product of
these specialization morphisms. We hope to return to this question
elsewhere. 

After that, we observe that there exists in the double Yangian
$DY(\SL_{2})$ a central element of the form 
$q^{(K+2)D}S$, deforming the difference $(K+2)D -
L_{-1}$, where $K$ and $D$ are the central and derivation elements of 
an extension of the loop algebra, $L_{-1}$ is a mode of the Sugawara
tensor, and $S$ belongs to the subalgebra of $DY(\SL_{2})$ ``without
$D$''. The construction of this element follows from a general
construction of V. Drinfeld in \cite{Dr} of central elements in
quasi-triangular Hopf algebras, implementing isomorphisms of modules
with their double duals. 

We then show that $S$ plays a role similar to $L_{-1}$ in the classical
situation, with infinitesimal shifts replaced by finite shifts of the
points: its copy $S^{(i)}$ on the $i$th factor of $DY(\SL_{2})^{\otimes n}$
conjugates the subalgebra $U_{\hbar}\G_{\zz}$ to $U_{\hbar}\G_{\zz + 
\hbar(K+2)\de_{i}}$, where $\de_{i}$ in the $i$th basis vector of
$\CC^{n}$ (see sect. \ref{act:Sug}).

As in the classical case, the actions of the $S^{(i)}$ define a
discrete flat
connection on the space of coinvariants $H_{0}(U_{\hbar}\G_{\zz},
\VV)^{*}$, where $\VV$ is the representation induced to
$DY(\SL_{2})^{\otimes n}$ from a finite-dimensional representation of
$U_{\hbar}\G_{\cO}$. We then compute explicitly this connection
(sect. \ref{cnx}), and find that it agrees with the quantum KZ
connection of \cite{FR}.

\section{Algebras associated with $\SL_{2}$ in the rational case}

\subsection{The algebra $U_{\hbar}\G_{\cK,\zz}$} \label{quant}

\subsubsection{Manin triples}  

Let us fix an integer $n\ge 1$. Let $z_{i}$, $i=1,\ldots,n$ be a family of
complex numbers; set $\zz = (z_{i})$. Let $t$ be the standard
coordinate on $\CC P^{1}$, and set $t_{i} = t-z_{i}$; $t_{i}$ is then a
local coordinate at $z_{i}$. 

Let $\cK_{i} = \CC((t_{i}))$ and $\cO_{i} = \CC[[t_{i}]]$ be the local
field and ring at $z_{i}$. Let us set 
$$
\cK = \oplus_{i=1}^{n} \cK_{i}, \quad \cO = \oplus_{i=1}^{n}\cO_{i}, 
$$
and let $R_{\zz}$ be the subring of $\cK$, formed by the expansions of
regular functions on $\CC P^{1} - \{z_{i}\}$, vanishing at infinity. 
A basis of $R_{\zz}$ is formed by the $(t-z_{i})^{-k-1}$, $i=
1,\ldots,n$, $k\ge 0$; as an element of $\cK$, $(t-z_{i})^{-k-1}$ should
be viewed as $((t_{j}+z_{j}-z_{i})^{-k-1})_{j=1,\ldots,n}$ . 
We then have the direct sum decomposition $\cK = R_{\zz} \oplus \cO$. 
Let us endow $\cK$ with the scalar product $\langle \phi,\psi \rangle_{\cK} =
\sum_{i=1}^{n}\res_{z_{i}}(\phi\psi dz)$. 

Let us set $\bar\G = \SL_{2}$, and construct the Lie algebras
$$
\G_{\cK,\zz} = (\bar\G \otimes \cK) \oplus \CC D_{\zz} \oplus \CC K_{\zz}
$$
as the double extension of the loop algebra $\bar\G\otimes \cK$ by the
cocycle $c(x\otimes \phi,y\otimes \psi) = \langle x, y\rangle_{\bar\G}
\sum_{i=1}^{n}\res_{z_{i}}(\phi d\psi)$, and by the derivation
$[D_{\zz},x\otimes \phi]
= x \otimes (d\psi/dz)$. 

This Lie algebra is endowed with the scalar product
$\langle , \rangle_{\G_{\cK,\zz}}$,
defined  on $\bar\G\otimes \cK$ by 
$$
\langle x\otimes\phi, y\otimes \psi\rangle_{\G_{\cK,\zz} } = \langle x,
y\rangle_{\bar\G}\langle\phi,\psi \rangle_{\cK}, 
$$
$\langle , \rangle_{\bar\G}$ being the Killing form of $\bar\G$, 
and $\langle D_{\zz},K_{\zz} \rangle_{\G_{\cK,\zz}} = 1$,  
$\langle D_{\zz},\bar\G\otimes \cK \rangle_{\wt\G_{\cK,\zz}} =
\langle K_{\zz},\bar\G\otimes \cK \rangle_{\wt\G_{\cK,\zz}} = 
\langle D_{\zz},D_{\zz} \rangle_{\G_{\cK,\zz}} 
\langle K_{\zz},K_{\zz} \rangle_{\G_{\cK,\zz}} 
=0$. 

The Lie algebra $\G_{\cK,\zz}$ contains subalgebras 
$$
\G_{\zz} = (\bar\G \otimes R_{\zz}) \oplus \CC K_{\zz} 
, \quad 
\G_{\cO} = (\bar\G \otimes\cO) \oplus \CC D_{\zz} . 
$$
Both subalgebras are isotropic for the scalar product $\langle\ ,
\ \rangle_{\G_{\cK,\zz}}$. 

We then construct the Manin triple 
\begin{equation} \label{triple}
\G_{\cK,\zz} = \G_{\zz} \oplus \G_{\cO}. 
\end{equation} 

In \cite{ER1}, we also considered the following twisted Manin
triples. Let $\bar\G = \N_{+} \oplus \HH \oplus \N_{-}$ be a Cartan
decomposition of $\bar\G$. Set
$$
\G_{+,\zz} =(\HH \otimes R_{\zz}) \oplus (\N_{+} \otimes \cK) \oplus 
\CC K_\zz, 
\quad
\G_{-} =(\HH \otimes \cO) \oplus (\N_{-} \otimes \cK) \oplus
\CC D_\zz , 
$$
and let $\G_{+,\zz}^{w_{0}} , \G_{\cO}^{w_{0}}$ be the subspaces
defined 
as $\G_{+,\zz} , \G_{\cO}$, exchanging $\N_{+}$ and
$\N_{-}$. Then we have the Manin triples
\begin{equation} \label{MT}
\G_{\cK,\zz} = \G_{+,\zz} 
\oplus \G_{-} , 
\end{equation}
and 
\begin{equation} \label{MTw0}
\G_{\cK,\zz}
 = \G_{+,\zz}^{w_{0}} 
\oplus \G_{-}^{w_{0}}. 
\end{equation}

\subsubsection{Quantization of the Manin triples (\ref{MT}) and
(\ref{MTw0})} 

In \cite{ER1}, we defined quantizations of the  Manin triples
(\ref{MT}) and $(\ref{MTw0})$. Let us recall their construction. 

Let $U_{\hbar}\G_{\cK,\zz}$ be the algebra generated by 
$D_{\zz},K_{\zz},x^{(i)}_{k}$, 
                                  $k\in\ZZ, i = 1,\cdots,n$, $x = e,f,h$,
arranged in the generating series 
$$
x_{+}(t) = \sum_{i=1}^n\sum_{k\ge 0} x^{(i)}_{k}(t-z_{i})^{-k-1}, 
\quad
x^{(i)}_{-}(t_{i}) = \sum_{k\ge 0} x^{[i]}_{-k-1}t_{i}^{k},  
$$
where
$$
x_{-k}^{[i]} = 
x_{-k}^{(i)} + \sum_{j\neq i,l\ge 0} (-1)^{l} {{l+k-1} \choose {k-1}}
z_{ji}^{-l-k}x_{l}^{(j)}
$$
(we note as usual $z_{ij} = z_{i} - z_{j}$); 
we also set 
$$
x^{(i)}(t_{i}) = \sum_{k\in\ZZ} x^{(i)}_{k} t_{i}^{-k-1};  
$$
we have then $\sum_{i=1}^{n} x^{(i)}(t_{i}) = 
x_{+}(z_{i} + t_{i}) + x^{(i)}_{-}(t_{i})$;
and we set  
\begin{equation} \label{K}
k_{+}(t) = \exp\left( {{q^{\pa} - 1}\over{\pa}} h_{+}(t)\right), \quad
K^{(i)}_{-}(t_{i}) = \exp \left( \hbar h^{(i)}_{-}(t_{i}) \right) ,
\quad \pa = d/dt, 
\end{equation}
and 
$K_{+}(t) = k_{+}(t_{i})k_{+}(t_{i}-\hbar)$, $K_{-}(t_{i}) =
k^{(i)}_{-}(t_{i})k^{(i)}_{-}(t_{i}-\hbar)$. 
In (\ref{K}), the arguments of the
exponentials are viewed as formal power series in $\hbar$, with
coefficients in 
$U_{\hbar}\G_{\cK,\zz}
\bar\otimes R_{\zz}$ 
in the first case, and in 
$U_{\hbar}\G_{\cK,\zz}\bar\otimes \cO$ in the second one.
Here $\bar\otimes$ denotes the graded tensor product with respect to the
bases $(t_{i}^{k})$ of $\cO$, and
$(t-z_{1})^{-k_{1}-1}\cdots(t-z_{n})^{-k_{n}-1}$ of $R_{\zz}$. 
In the notation of \cite{ER1}, we have $x^{(i)}_{k} = x[t_{i}^{k}]$, 
$x^{[i]}_{k} = x[(t-z_{i})^{k}]$. 
Thus $x^{(i)}_k$ is a $q$-analogue of $x\otimes t_i^k$ in
$\bar\G\otimes\cK_i\subset\bar\G\otimes\cK$ and $x^{[i]}_k$ is
a $q$-analogue of $x\otimes(t-z_i)^k\in \bar\G\otimes R_{\zz}$.
The above relation between $x^{[i]}_k$ and the generators $x^{(j)}_l$
is obtained by the Laurent expansion at $z_j$.

The relations are 
\begin{equation} \label{k+:e}
(t-z_{j} - u_{j}) k_{+}(t)e^{(j)}(u_{j}) = 
(t-z_{j} - u_{j}+\hbar) e^{(j)}(u_{j}) k_{+}(t) ,
\end{equation}
\begin{equation} \label{k+:f}
(t-z_{j} - u_{j}+\hbar)
k_{+}(t)f^{(j)}(u_{j})
= 
(t-z_{j} - u_{j})
f^{(j)}(u_{j}) k_{+}(t),
\end{equation}
\begin{equation} \label{k-:e}
(z_{ij} + t_{i} - u_{j} - \hbar K_{\zz} + \hbar) 
k^{(i)}_{-}(t_{i})e^{(j)}(u_{j})
= (z_{ij} + t_{i} - u_{j}-\hbar K_{\zz})
e^{(j)}(u_{j}) k^{(i)}_{-}(t_{i}) ,
\end{equation}
\begin{equation} \label{k-:f}
(z_{ij} + t_{i} - u_{j}) k^{(i)}_{-}(t_{i})f^{(j)}(u_{j})
= (z_{ij} + t_{i} - u_{j}+\hbar) f^{(j)}(u_{j}) k^{(i)}_{-}(t_{i}),
\end{equation}
\begin{equation} \label{e:e}
(z_{ij} + t_{i} - u_{j}- \hbar)e^{(i)}(t_{i})e^{(j)}(u_{j}) 
= (z_{ij} + t_{i} - u_{j}+ \hbar) e^{(j)}(u_{j}) e^{(i)}(t_{i}), 
\end{equation}
\begin{equation} \label{f:f}
(z_{ij} + t_{i} - u_{j}+ \hbar)f^{(i)}(t_{i})f^{(j)}(u_{j}) 
= (z_{ij} + t_{i} - u_{j}- \hbar) f^{(j)}(u_{j}) f^{(i)}(t_{i}), 
\end{equation}
\begin{equation} \label{e:f}
[e^{(i)}(t_{i}),f^{(j)}(u_{j})] = 
\de_{ij}{1\over{\hbar}} \left(\delta(t_{i},u_{j})
K^{+}(z_{i}+t_{i}) - 
\delta(t_{i},u_{j}-\hbar K_{\zz})
K^{(j)}_{-}(u_{j})^{-1} \right) , 
\end{equation}
\begin{equation} \label{k:k}
[K_{\zz},\on{anything}] = [k_{+}(t), k_{+}(u)] = [K_{-}^{(i)}(t_{i}) , 
K_{-}^{(j)}(t_{j})] =  0, 
\end{equation}
\begin{align} \label{K+:K-}
(t - z_{j} - & u_{j}-\hbar)(t-z_{j} - u_{j}+\hbar K_{\zz} + \hbar )
 K_{+}(t) K^{(j)}_{-}(u_{j}) 
\\ & \nonumber = (t - z_{j} - u_{j}+\hbar)
(t - z_{j} - u_{j}+\hbar K_{\zz} - \hbar)
K^{(j)}_{-}(u_{j}) K_{+}(t).  
\end{align}
\begin{equation} \label{D:any}
[D_{\zz},x^{(i)}(t_{i})] = - ({dx^{(i)}}/{dt_{i}})(t_{i}), \quad x=e,f,h. 
\end{equation}

Here $\de(z,w) =\sum_{k\in\ZZ} z^{k}w^{-k-1}$; note 
that we have changed both signs of 
$K_{\zz}$ and $D_{\zz}$ with respect to 
the convention of \cite{EF}.

This algebra is endowed with Hopf structures $(\Delta,\varepsilon,S)$ 
and $(\bar \Delta,\varepsilon,\bar S)$, quantizing the bialgebra
structures (\ref{MT}) and (\ref{MTw0}) respectively.
These coproducts are defined by   
\begin{equation} \label{Delta:k}
\Delta(k_{+}(t)) = k_{+}(t) \otimes k_{+}(t), \quad
\Delta(K^{(i)}_{-}(t_{i})) = K^{(i)}_{-}(t_{i}) \otimes K^{(i)}_{-}(t_{i} -
\hbar (K_{\zz})_{1}), 
\end{equation}
\begin{equation} \label{Delta:e}
\Delta(e^{(i)}(t_{i})) = e^{(i)}(t_{i})\otimes K_{+}(z_{i}+t_{i}) +
1\otimes e^{(i)}(t_{i}), 
\end{equation}
\begin{equation} \label{Delta:f}
\Delta(f^{(i)}(t_{i})) = f^{(i)}(t_{i})\otimes 1 + K_{-}^{(i)}(t_{i})^{-1} 
\otimes f^{(i)}(t_{i}-\hbar (K_{\zz})_{1}),
\end{equation}
\begin{equation} \label{Delta:D:K}
\Delta(D_{\zz}) = D_{\zz} \otimes 1 + 1\otimes D_{\zz}, 
\quad \Delta(K_{\zz}) =  K_{\zz}\otimes 1 + 1 \otimes K_{\zz}, 
\end{equation}
where $(K_{\zz})_{1}$ and $(K_{\zz})_{2}$ mean $K_{\zz}\otimes 1$ and $1
\otimes K_{\zz}$; and  
\begin{equation} \label{bar:Delta:k}
\bar\Delta(k_{+}(t)) = k_{+}(t) \otimes k_{+}(t), 
\quad
\bar\Delta(K^{(i)}_{-}(t_{i})) = K^{(i)}_{-}(t_{i})
\otimes K^{(i)}_{-}(t_{i} + \hbar (K_{\zz})_{1})), 
\end{equation}
\begin{equation} \label{bar:Delta:e}
\bar\Delta(e^{(i)}(t_{i})) = e^{(i)}(t_{i} - \hbar (K_{\zz})_{2})\otimes
K^{(i)}_{-}(t_{i}- \hbar (K_{\zz})_{2})^{-1} 
+ 1\otimes e^{(i)}(t_{i}),
\end{equation}
\begin{equation} \label{bar:Delta:f}
\bar\Delta(f^{(i)}(t_{i})) = f^{(i)}(t_{i})\otimes 1 + K^{(i)}_{+}(t_{i}) 
\otimes f^{(i)}(t_{i}), 
\end{equation}
\begin{equation} \label{bar:Delta:D:K}
\bar\Delta(D_{\zz}) = D_{\zz} \otimes 1 + 1\otimes D_{\zz}, \quad
\bar\Delta(K_{\zz}) = K_{\zz} \otimes 1 + 1 \otimes K_{\zz}.  
\end{equation}

The counit $\varepsilon$ is defined to be equal to zero on all
generators. 

\begin{prop} (see \cite{ER1})
The above formulas define quantizations $(U_{\hbar}\G_{\cK,\zz},
\Delta)$ and $(U_{\hbar}\G_{\cK,\zz},\bar\Delta)$ 
of the Manin triples (\ref{MT}) and (\ref{MTw0}). 
These are quasitriangular Hopf algebras, with universal
$R$-matrices respectively given by 
$$
\cR = q^{D_{\zz} \otimes K_{\zz}}\exp\left( {\hbar\over 2} 
\sum_{i=1}^{n}\sum_{k\ge 0} 
h^{(i)}_{k} \otimes h^{[i]}_{-k-1}
\right) \exp \left( \hbar \sum_{i=1}^{n}\sum_{k\in\ZZ} 
e^{(i)}_{k} \otimes f^{(i)}_{-k-1}
 \right) , 
$$
$$
\bar\cR = \exp \left( \hbar \sum_{i=1}^{n}\sum_{k\in\ZZ}
f^{(i)}_{k} \otimes e^{(i)}_{-k-1}  \right)
 q^{D_{\zz} \otimes K_{\zz}}
\exp\left( {\hbar\over 2} \sum_{i=1}^{n}\sum_{k\ge 0} 
h^{(i)}_{k} \otimes h^{[i]}_{-k-1}
\right) . 
$$
\end{prop}

The coproducts $\Delta$ and $\bar\Delta$ are related by the twist 
$$
F = \exp\left( \hbar \sum_{i=1}^{n} \sum_{k\in\ZZ} e^{(i)}_{k} \otimes
f^{(i)}_{-k-1}\right);  
$$
this means that we have 
$$
\bar\Delta = \Ad(F) \circ \Delta. 
$$
We denote $\Ad(X)$, for an invertible element $X$ of some algebra $A$,
the linear endomorphism $Y\mapsto XYX^{-1}$ of $A$.
$F$ also satisfies the cocycle equation 
\begin{equation} \label{cocycle}
(F\otimes 1)(\Delta\otimes 1)(F)= 
(1 \otimes F)(1\otimes \Delta)(F).
\end{equation}

\begin{remark} \label{extension}
The above results can be extended to the case where we take 
$\zz$ in 
$\CC[[\hbar]]^{n}$, such that $z_{i}-z_{j}\notin \hbar\CC[[\hbar]]$ for
$i\neq j$. \hfill \qed\medskip 
\end{remark} 

\subsubsection{Subalgebras of $U_{\hbar}\G_{\cK,\zz}$}

Let us consider in $U_{\hbar}\G_{\cK,\zz}$, the subalgebras
$U_{\hbar}\G_{\cO}$ and $U_{\hbar}\G_{\zz}$, respectively generated by
$D_{\zz}$ and the $x^{(i)}_{k}$, $k\ge 0$, and by $K_{\zz}$ and the 
$$
x_{-k}^{[i]} = 
x_{-k}^{(i)} + \sum_{j\neq i,l\ge 0} (-1)^{l} {{l+k-1} \choose {k-1}}
z_{ji}^{-l-k}x_{l}^{(j)}, \quad k\ge 1;  
$$
$x=e,f,h$.

\begin{prop} (see \cite{ER2}, sect. 4)
The inclusion of algebras $U_{\hbar}\G_{\zz} \subset 
U_{\hbar}\G_{\cK,\zz}$  and $U_{\hbar}\G_{\cO} \subset 
U_{\hbar}\G_{\cK,\zz}$ 
are flat deformations of the inclusions of enveloping
algebras $U\G_{\zz} \subset U\G_{\cK,\zz}$ and
$U\G_{\cO} \subset U\G_{\cK,\zz}$. 
\end{prop}

\begin{prop} \label{PBW}
The product maps 
$$
U_{\hbar}\G_{\zz} \otimes U_{\hbar}\G_{\cO} \to U_{\hbar} \G_{\cK,\zz}
\ \on{and}\ 
U_{\hbar}\G_{\cO}  \otimes U_{\hbar}\G_{\zz}
\to U_{\hbar} \G_{\cK,\zz}
$$
define linear isomorphisms.  
\end{prop} 

{\em Proof.} This follows from the fact that the three spaces
$U_{\hbar}\G_{\zz}, U_{\hbar}\G_{\cO}$ and $U_{\hbar} \G_{\cK,\zz}$ 
are flat deformations of $U\G_{\zz}, U\G_{\cO}$ and $U\G_{\cK,\zz}$, and that 
the product maps $U\G_{\zz} \otimes U\G_{\cO} \to U\G_{\cK,\zz}$ and
$U\G_{\cO}  \otimes U\G_{\zz} \to U\G_{\cK,\zz}$ 
are linear isomorphisms. 
\hfill \qed

\subsection{Hopf structures on $U_{\hbar}\G_{\cK,\zz}$ and its subalgebras} 

In a way similar to \cite{EF}, we define the linear maps
$$
\Pi_{\zz,\ell}: U_{\hbar}\G_{\cK,\zz} \to U_{\hbar}\G_{\zz}, \quad
\Pi_{\cO,r}: U_{\hbar}\G_{\cK,\zz} \to U_{\hbar}\G_{\cO}
$$
by posing 
$$
\Pi_{\zz,\ell}(a_{\zz}a_{\cO}) = a_{\zz}\varepsilon(a_{\cO}), \quad
\Pi_{\cO,r}(a_{\zz}a_{\cO}) = \varepsilon( a_{\zz}) a_{\cO}, 
$$
and 
$$
\Pi_{\cO,\ell}: U_{\hbar}\G_{\cK,\zz} \to U_{\hbar}\G_{\cO}, \quad
\Pi_{\zz,r}: U_{\hbar}\G_{\cK,\zz} \to U_{\hbar}\G_{\zz}
$$
by posing 
$$
\Pi_{\cO,\ell}(a_{\cO}a_{\zz}) = a_{\cO}\varepsilon(a_{\zz}), \quad
\Pi_{\zz,r}(a_{\cO}a_{\zz}) = \varepsilon( a_{\cO}) a_{\zz}, 
$$
for $a_{\zz}\in U_{\hbar}\G_{\zz}, a_{\cO}\in U_{\hbar}\G_{\cO}$. The fact 
that these maps are well-defined follows from Prop. \ref{PBW}. 

We then have 

\begin{thm} (see \cite{ER1,EK}) We have the equalities 
$$
(\Pi_{\zz,r} \otimes 1)(F) = (1 \otimes \Pi_{\cO,r})(F) \ \on{and} \ 
(\Pi_{\cO,\ell}\otimes 1)(F) = (1\otimes \Pi_{\zz,\ell})(F). 
$$
If we set 
$$
F_{1} = (\Pi_{\zz,r} \otimes 1)(F) = (1 \otimes \Pi_{\cO,r})(F) , 
$$ 
and 
$$
F_{2} = (\Pi_{\cO,\ell}\otimes 1)(F) = (1\otimes \Pi_{\zz,\ell})(F), 
$$
we have $F = F_{2}F_{1}$. The maps $\Ad(F_{1}) \circ\Delta$ and
$\Ad(F_{2}^{-1})\circ\bar\Delta$ coincide; we will denote them by 
$\Delta_{\zz}$. $\Delta_{\zz}$ defines a quasitriangular
Hopf algebra structure on $U_{\hbar}\G_{\cK,\zz}$; 
$U_{\hbar}\G_{\zz}$ and $U_{\hbar}\G_{\cO}$ are Hopf subalgebras of 
$U_{\hbar}\G_{\cK,\zz}$ for this structure. The universal $R$-matrix of 
$(U_{\hbar}\G_{\cK,\zz} , \Delta_{\zz})$ is expressed by 
\begin{equation} \label{R:mat}
\cR_{\zz} = F_{1}^{(21)} q^{D_{\zz} \otimes K_{\zz}} 
\exp\left( {\hbar\over 2} \sum_{i=1}^{n}\sum_{k\ge 0} 
h^{(i)}_{k} \otimes h^{[i]}_{-k-1}
\right) F_{2}
\end{equation}
\end{thm}

Moreover, $F_1$ and $F_2$ have the expansions 
\begin{equation} \label{expansion}
F_1 \in 1 + \hbar \sum_{i=1}^n\sum_{k\ge 0}e^{(i)}_{k} \otimes 
f^{[i]}_{-k-1}
+ U_{\hbar}\N_{+}^{[2]} \otimes U_{\hbar}\N_{-}^{[2]}
\end{equation}
and 
\begin{equation} \label{expansion2}
F_2 \in 1 + \hbar \sum_{i=1}^n\sum_{k\ge 0}e^{[i]}_{-k-1} \otimes 
f^{(i)}_{k}
+ U_{\hbar}\N_{+}^{[2]} \otimes U_{\hbar}\N_{-}^{[2]}, 
\end{equation}
where $U_{\hbar}\N_{\pm}^{[2]}$ are the linear spans in 
$U_{\hbar}\G_{\cK,\zz}$ of products of 
more than two elements of the form $x[\eps]$,  
$x=e$ for $\pm = +$ and $x=f$ for $\pm = -$, $\eps\in\cK$.  

The Hopf algebra $(U_{\hbar}\G_{\cK,\zz}, \Delta_{\zz})$, together
with its Hopf subalgebras $U_{\hbar}\G_{\zz}$ and $U_{\hbar}\G_{\cO}$, 
forms a quantization of the Manin triple (\ref{triple}).

\subsection{Representations of $U_{\hbar}\G_{\cK,\zz}$ and $L$-operators}

In \cite{ER1}, we studied level zero representations of the algebras
introduced there. Our result can be expressed as follows. Define the
   algebra $U_{\hbar}\G'_{\cK,\zz}$ as the
                                                       subalgebra 
of $U_{\hbar}\G_{\cK,\zz}$ generated by the $K_{\zz}$ and the 
$x[\eps],\eps\in \cK,x=e,f,h$; and 
   $U_{\hbar}\G_{\zz}$ and $U_{\hbar}\G'_{\cO}$
as its subalgebras respectively generated by the $x[\eps],\eps\in
R_{\zz}$ 
and $K_\zz$,
 and by the $x[\eps],\eps\in\cO$, $x = e,f,h$. 

\begin{prop} \label{repr}
We have an algebra morphism 
$$
\pi : U_{\hbar}\G'_{\cK,\zz} \to
\End(\CC^{2})  \otimes \cK[[\hbar]], 
$$
defined by
$$
\pi (h_{k}^{(i)}) =   
\pmatrix  {2\over{1+q^{\pa}}}t_{i}^{k} 
 & 0 \\ 0 & 
-  {2\over{1+q^{-\pa}}}t_{i}^{k} \endpmatrix, \quad
k\ge 0,
$$
$$
\pi (h^{[i]}_{k}) = \pmatrix 
\left( {{1-q^{-\pa}}\over{\hbar\pa}} \right) (t-z_{i})^{k}
& 0 \\ 0 & 
- \left( {{q^{\pa}-1}\over{\hbar\pa}} (t-z_{i})^{k}\right) 
\endpmatrix, \quad k <0,
$$
$$
\pi (e^{(i)}_{k}) = 
\pmatrix 0 & t_{i}^{k} \\ 0 & 0 \endpmatrix, 
\quad
\pi (f^{(i)}_{k}) = 
\pmatrix 0 & 0 \\ t_{i}^{k} & 0 \endpmatrix, \quad
k\in\ZZ, 
$$
and $\pi(K_\zz)=0$,
where $\pa$ is the derivation of $\cK$ defined as 
$\sum_{i=1}^n d/dt_i$; 
it coincides with $d/dt$ when restricted to $R_\zz$. 

The images by $\pi$ of $U_{\hbar}\G_{\zz}$ and
$U_{\hbar}\G'_{\cO}$ 
are contained in $\Id_{\CC^2}\otimes\CC[[\hbar]]+
\hbar\End(\CC^{2}) \otimes R_{\zz}[[\hbar]]$ 
and  $\End(\CC^{2}) \otimes \cO[[\hbar]]$  respectively. 
 
\end{prop}

Define the $L$-operators of $(U_\hbar\G_{\cK,\zz} , \Delta_\zz)$ 
as 
$$
L^{+}_\zz= (\pi\otimes 1)(\cR_{\zz}^{(21)}), \quad
L^{-}_\zz = (\pi\otimes 1)(\cR_{\zz}^{-1} q^{D_{\zz} \otimes K_{\zz}}).
$$
Since $\cR_\zz$ is contained in the completion of $U_\hbar\G_\cO\otimes
U_\hbar\G_\zz$, see (\ref{R:mat}), and has leading term $1$,
we have 
$$
L^{+}_\zz \in 1+\hbar(\End(\CC^{2})\otimes R_{\zz})\bar\otimes U_{\hbar}
                                                                    \G_{\cO},
\quad 
L^{-}_\zz \in 1+ \hbar(\End(\CC^{2})\otimes \cO)\bar\otimes U_{\hbar}
                                                                   \G_{\zz}.
$$
In what follows, we will stress the functional dependences 
of $L^{+}_\zz$ and $L^{-}_\zz$ by writing $L^{+}_\zz$ as
$L^{+}_\zz(t)$ and the component of  $L^{-}_\zz$ in $\cO_i$
as $L^{-}_{(i)}(t_i)$.

{}From (\ref{expansion}) and (\ref{expansion2}) 
follows that $L^{+}_\zz(t)$ and $L^{-}_{(i)}(t_i)$ 
are decomposed as follows: 
\begin{equation} \label{L+}
L^{+}_\zz(t) = 
\pmatrix 1 & \hbar f_{+}(t) \\ 0 & 1 \endpmatrix 
\pmatrix k_{+}(t-\hbar) & 0 \\ 0 & k_{+}^{-1}(t) \endpmatrix
\pmatrix 1 & 0 \\ \hbar e_{+}(t) & 1 \endpmatrix
, 
\end{equation}
and
\begin{equation} \label{L-}
L^{-}_{(i)}(t_{i}) = 
\pmatrix 1 & - \hbar f^{(i)}_{-}(t_{i}) \\ 0 & 1\endpmatrix 
\pmatrix k_{-}^{(i)}(t_{i}-\hbar)^{-1} & 0 \\ 0 & k_{-}^{(i)}(t_{i}) 
\endpmatrix
\pmatrix 1 & 0 \\ - \hbar e^{(i)}_{-}(t_{i}+\hbar K_{\zz}) & 1
\endpmatrix . 
\end{equation}

\subsection{$RLL$ relations for $U_{\hbar}\G_{\cK,\zz}$ and its
subalgebras} \label{rll}

Let us compute the image by $\pi\otimes\pi$ of the universal $R$-matrix
of $U_{\hbar}\G_{\cK,\zz}$. 

We find
\begin{equation} \label{im:R}
(\pi\otimes \pi)(\cR_{\zz})(t_{i},u) = R^{-}(z_i + t_i , u) , 
\end{equation}
where 
\begin{align} \label{R-}
R^{-}(z_i + t_i , u) & = \exp\left( \sum_{k\ge 0} ( {1\over\pa}{ {q^\pa -
1} \over  
{ q^\pa + 1}} t_i^k) (u-z_i)^{-k-1}\right) \cdot 
\\ & \nonumber \cdot
{ 1\over{z_i +t_i -u -\hbar}}(
(z_i +t_i -u)\Id_{\CC^2 \otimes \CC^2} - \hbar P) ,  
\end{align}  
expanded at the vicinity of $t_i = 0$, 
where $P$ is the permutation operator of the two 
factors of $\CC^2 \otimes \CC^2$. 

By applying the representation $\pi$ to two out of the three factors of 
the Yang-Baxter equation 
$$
\cR_{\zz}^{(12)} \cR_{\zz}^{(13)} \cR_{\zz}^{(23)}
= \cR_{\zz}^{(23)}\cR_{\zz}^{(13)}\cR_{\zz}^{(12)}, 
$$
we obtain the relations
\begin{equation} \label{RLL++:zz}
R(t,t')L^{+(1)}_{\zz}(t)L^{+(2)}_{\zz}(t')= 
L^{+(2)}_{\zz}(t') L^{+(1)}_{\zz}(t) R(t,t') , 
\end{equation}
\begin{align} \label{RLL+-:zz}
R(z_{i} + t_{i} , t' )
L^{-(1)}_{(i)}(t_{i}) 
L^{+(2)}_{\zz}(t')
= L^{+(2)}_{\zz}(t') L^{-(1)}_{(i)}(t_{i}) 
R(z_{i} + t_{i} , t'+\hbar K_{\zz}) , 
\end{align}
\begin{equation} \label{RLL--:zz}
R(z_{ij} + t_{i} , t_{j}) 
L^{-(1)}_{(i)}(t_{i}) 
L^{-(2)}_{(j)}(t_{j}) 
= L^{-(2)}_{(j)}(t_{j}) L^{-(1)}_{(i)}(t_{i})  
R(z_{ij} + t_{i} , t_{j}) , 
\end{equation}
with 
\begin{equation} \label{R}
R(z,z') = { 1\over{z-z'+\hbar}}
( (z-z') \Id_{\CC^2 \otimes \CC^2} + \hbar P) \exp\left( -\sum_{k\ge 0}
\left( ({1\over \pa}{{q^{\pa}-1}\over{q^{\pa}+1}})z^{k} \right)  
z^{{\prime
-k-1}}\right) ;  
\end{equation}
we set 
$$
R(z) = { 1\over{z+\hbar}}
( z\Id_{\CC^2 \otimes \CC^2} + \hbar P) \exp\left( 
({1\over \pa}{{q^{\pa}-1}\over{q^{\pa}+1}})z^{-1} \right), 
$$
so that 
           $R(z,z') = R(z-z')$. 

Moreover, we have the relations
\begin{equation} \label{det}
\det_{\hbar}(L^{+}_{\zz}(t)) = \det_{\hbar}(L^{-}_{(i)}(t_{i})) = 1,  
\end{equation}
where 
$$
\det_{\hbar}(L(t)) = l_{11}(t)l_{22}(t-\hbar) - l_{12}(t)l_{21}(t-\hbar),
$$
for $L(t)$ a matrix with entries $l_{ij}(t)$; this follows from the 
identities 
$$
k^{+}(t) e^{+}(u) k^{+}(t)^{-1} = {{t-u+\hbar}\over{t-u}}e^{+}(u)
- {{\hbar}\over{t-u}}e^{+}(t), 
$$
$$
k^{-}(t_{i})^{-1} e^{-}(u_{i}+\hbar K_{\zz}) k^{-}(t_{i})^{-1} = 
{{t_{i}-u_{i}+\hbar}\over{t_{i}-u_{i}}}e^{-}(u_{i}-\hbar K_{\zz})
- {{\hbar}\over{t_{i}-u_{i}}}e^{-}(t_{i}),
$$
specialized to $t = u-\hbar$ and $t_{i} = u_{i}-\hbar$.

Below we will denote by $1_{i}$ the element of $\cK$ with $i$th
component equal to $1$ and the others equal to $0$. 

\begin{prop}
The algebra with generators $D_{\zz}, K_{\zz}$ and 
$l^{(i)}_{\al\beta}[k]$, $\al,\beta = 1,2, 
k\in\ZZ, i = 1, \ldots ,n$, 
arranged in matrices 
$$
L^{+}_{\al\beta}(t) = 
\de_{\al\beta} + \hbar \sum_{i=1}^{n} \sum_{k\ge 0} l^{(i)}_{\al\beta}[k]
(t-z_{i})^{-k-1}
, \quad
L^{-(i)}_{\al\beta}(t_{i}) = \de_{\al\beta}1_{i} + 
\hbar \sum_{k\ge 0} l^{(i)}_{\al\beta}[-k-1] t_{i}^{k} ,  
$$
subject to relations (\ref{RLL++:zz}), (\ref{RLL+-:zz}), 
(\ref{RLL--:zz}), (\ref{det}), and  
\begin{equation} \label{rest}
[K_{\zz}, \on{anything}] = 0, \quad [D_{\zz}, L^{+}_{\al\beta}(t) ]
= - { { dL^{+}_{\al\beta}(t)} \over{dt}}, \quad  
[D_{\zz}, L^{-(i)}_{\al\beta}(t_{i}) ]
= - { { dL^{-(i)}_{\al\beta}(t_{i})} \over{dt_{i}}}, 
\end{equation}
is isomorphic to 
$U_{\hbar}\G_{\cK,\zz}$. 
\end{prop}

{\em Proof.} Denote by $U$ the algebra defined above. 
As we have seen before, the formulas
(\ref{L+}) and (\ref{L-}) define a morphism from $U$ 
to $U_{\hbar}\G_{\cK,\zz}$. 
On the other hand, because of relations (\ref{det}), a system of 
generators of $U$ is given by $D_{\zz}, K_{\zz}$ and
the $l_{\al\beta}^{(i)}[k]$, $(\al,\beta)
\neq (2,2)$, $i=1,\ldots,n, k\in \ZZ$. 
Relations (\ref{RLL++:zz}), (\ref{RLL+-:zz})
and (\ref{RLL--:zz}) allow to write commutators 
$ [ l_{\al\beta}^{(i)}[k] , l_{\al'\beta'}^{(i')}[k'] ]$ 
as linear combinations of the 
$l_{\al\beta}^{(i)}[k]$, up to first order in $\hbar$; it is easy to 
check that this combination is given by the Lie algebra bracket in 
$\bar\G \otimes \cK$. 
Therefore, $U$ is
a deformation of the tensor product of the symmetric
algebra of $\bar\G \otimes \cK$
with the symmetric algebra in $D_{\zz}$ and $K_{\zz}$. 
In the classical limit, the morphism defined by 
(\ref{L+}) and (\ref{L-}) is the identity. 
\hfill \qed\medskip 

\begin{corollary} \label{funct}
Let $\varphi$ be an algebra morphism from $U_{\hbar}\G_{\cK,\zz}$
to some algebra $\cA$. Then the images by $\varphi\otimes 1 \otimes 1$
of $L^{+}(t)$ and $L^{-}_{(i)}(t_{i})$ are elements
\begin{equation} \label{phi:L+}
\varphi(L^{+}(t)) \in 1\otimes 1 + \hbar (\cA \otimes \End(\CC^{2}))\bar
\otimes R_{\zz}
               [[\hbar]], 
\end{equation}
\begin{equation} \label{phi:L-}
\varphi(L^{-}_{(i)}(t_{i})) \in 1\otimes 1_{i} +\hbar 
(\cA \otimes \End(\CC^{2}))\bar\otimes \cO
                                          [[\hbar]], 
\end{equation} 
satisfying (\ref{RLL+-:zz}), (\ref{RLL+-:zz}), (\ref{det}). Conversely, 
any matrices (\ref{phi:L+}), (\ref{phi:L-}) 
 satisfying these equations define a
morphism $\varphi: U_{\hbar}\G_{\cK,\zz}\to \cA$.  
\end{corollary}

\subsection{Isomorphism of $U_{\hbar}\G_{\cK,\zz}$ with
$DY(\SL_{2})^{\otimes n}/(K^{(i)}-K^{(j)})$}

In this section, we will construct an isomorphism $i_{\zz}$ from  
$U_{\hbar}\G_{\cK,\zz}$ to 
$DY(\SL_{2})^{\otimes n}/(K^{(i)}-K^{(j)})$. We set 
$K^{(i)} = 1^{\otimes (i-1)} \otimes K \otimes 1^{\otimes (n-i)}$ 
and we will denote by $K$ the common value of the 
$K^{(i)}$ in the latter algebra. 

\subsubsection{Presentation of $DY(\SL_2)$}

We will express $i_{\zz}$ in terms of $L$-operators. To do that, 
let us remark that there is a simple presentation of $DY(\SL_{2})$
as $U_{\hbar}\G_{\CC((t)),(0)}$, the specialization for $n= 1$ and 
$\zz = (0)$ of the algebra $U_{\hbar}\G_{\cK,\zz}$. 

In terms of $L$-operators, $DY(\SL_{2})$ is presented as follows. 
It has generators $D,K$ and $l_{\al\beta}[k]$, $\al,\beta = 1,2, k\in\ZZ$, 
generating series 
$$
L^{+}_{\al\beta}(z) = \de_{\al\beta} + \hbar \sum_{k\ge 0} l_{\al\beta}[k]
z^{-k-1}, 
\quad 
L^{-}_{\al\beta}(z) = \de_{\al\beta} + \hbar \sum_{k\ge 0} 
l_{\al\beta}[-k-1] z^{k}, 
$$
and relations 
\begin{equation} \label{RLLpm:Y}
R(z , z') L^{\pm(1)}(z)L^{\pm(2)}(z')
= 
L^{\pm(2)}(z')L^{\pm(1)}(z)R(z , z'), 
\end{equation}
\begin{equation} \label{RLL+-:Y}
R(z , z') L^{-(1)}(z)L^{+(2)}(z')
= 
L^{+(2)}(z')L^{-(1)}(z)R(z , z'-\hbar K),
\end{equation}
where $R(z,z')$ is defined by 
(\ref{R}), 
the quantum determinant relations
\begin{equation} \label{det:Y}
\det_{\hbar}(L^{\pm}(z)) =  1,
\end{equation}
and 
$$
[K,\on{anything}] = 0, \quad 
[D,L^{\pm}(z)] = - dL^{\pm}(z)/dz. 
$$
(In the notation of \cite{EF}, 
$L^{-}(z)$ and $L^{+}(z)$ correspond respectively to $L^{<0}(z)$ 
and to the inverse of $L^{\ge 0}(z)$ in $\End(\CC^{2})[[z^{\pm 1}]]
\bar\otimes DY(\SL_{2})$.)

We also have the relation
$$
R(z,z') L^{+(1)}(z)L^{-(2)}(z')
= 
L^{-(2)}(z')L^{+(1)}(z)R(z,z'+\hbar K). 
$$

\subsubsection{The isomorphism}

Let $L^{\pm}_{i}(z)$ be the image in 
$$
DY(\SL_{2})^{\otimes n} \otimes
\End(\CC^{2})[[z^{\pm1}]][[\hbar]]
$$
of the operator $L^{\pm}(z)$, which
belongs to $DY(\SL_{2}) \otimes \End(\CC^{2})[[z^{\pm1}]][[\hbar]]$, by
the map $\ell \otimes a\mapsto 1^{\otimes (i-1)}\otimes \ell \otimes
1^{\otimes (n-i)} \otimes a$. 

Consider the expression 
$$
L^{+}_{1}(t_{i} + z_{i1} -\hbar K)\cdots
L^{+}_{i-1}(t_{i}+z_{i,i-1} -\hbar K)
L^{-}_{i}(t_{i})
L^{+}_{i+1} (t_{i}+z_{i,i+1}) \cdots 
L^{+}_{n}(t_{i}+z_{i,n}). 
$$
$L^{-}_{i}(t_{i})$ belongs to $1\otimes 1_{i} +\hbar 
DY(\SL_2)^{\otimes n}[[t_i]]$. On the other
hand, the series $L^{+}_{j}(t_{i} + z_{ij}),j\neq i$ 
are expanded as sums
$1\otimes 1 +\hbar\sum_{k\ge 1}(t_{i} + z_{ij})^{-k}\times\on{coefficients}$, 
and we expand them in turn at the
vicinity of $t_{i} = 0$. The coefficient of each power of $t_{i}$ in this
expansion is then an infinite series, converging in the topology of
$DY(\SL_{2})^{\otimes n}$ because it involves large Fourier indices. 
Therefore the above expression belongs to $1\otimes 1_{i} + \hbar
DY(\SL_2)^{\otimes n}[[t_{i}]]$. 

Consider now the expression 
$$
L^{+}_{1}(t - z_{1})\cdots L^{+}_{n}(t - z_{n}). 
$$
Since each $L^{+}(t)$ belongs to $1\otimes 1 + \hbar t^{-1}DY(\SL_2)
[[t^{-1}]]$, this product belongs to $1\otimes 1 + 
\hbar DY(\SL_2)^{\otimes n} \bar\otimes 
R_{\zz}$. 

Then 

\begin{prop} \label{isom}
The formulas
\begin{align} \label{isom:L-}
i_{\zz}( L^{-}_{(i)}(t_{i}))  = L^{+}_{1}(t_{i} + z_{i1})\cdots
L^{+}_{i-1}(t_{i}+z_{i,i-1})
L^{-}_{i}(t_{i})
L^{+}_{i+1} & (t_{i}+z_{i,i+1} - \hbar K) \cdots \nonumber \\ & 
L^{+}_{n}(t_{i}+z_{i,n} - \hbar K),
\end{align}
\begin{align} \label{isom:L+}
i_{\zz}( L^{+}_{\zz}(t)) = 
L^{+}_{1}(t - z_{1})\cdots L^{+}_{n}(t - z_{n}), 
\end{align}
\begin{align} \label{isom:rest}
i_{\zz}(D_{\zz}) =  \sum_{i=1}^{n}D^{(i)}, \quad i_{\zz}(K_{\zz}) =  K,  
\end{align}
define an isomorphism between $U_{\hbar}\G_{\cK,\zz}$ and 
$DY(\SL_2)^{\otimes n}/(K^{(i)} - K^{(j)})$. 
\end{prop}

{\em Proof. }
By the above remarks, the right sides of (\ref{isom:L-}) and 
(\ref{isom:L+}) have the form prescribed by Cor. \ref{funct}. 
We then check that the right sides of formulas (\ref{isom:L-})
and (\ref{isom:L+}) satisfy relations (\ref{RLL++:zz}), (\ref{RLL+-:zz})
and (\ref{RLL--:zz}). 
This follows from the Yangian relations 
(\ref{RLLpm:Y}) and (\ref{RLL+-:Y}). 

The fact that the right sides of (\ref{isom:L-}) and 
(\ref{isom:L+}) satisfy the identities (\ref{det}) follows from the
fact that the quantum determinant is a group-like element of the Yangian
algebra $Y(\GL_2)$. 

The relations (\ref{rest}) are obviously satisfied by the 
right sides of (\ref{isom:rest}), (\ref{isom:L-}) and
(\ref{isom:L+}). 
\hfill \qed\medskip

\subsection{Shifts of the points} \label{shifts}

By Prop. \ref{isom}, we have subalgebras $i_{\zz}(U_{\hbar}\G_{\zz})$
of 
$$
DY(\SL_{2})^{\otimes n}/(K^{(i)} - K^{(j)}),
$$ 
for $\zz \in
\CC[[\hbar]]^{n}$, such that $z_{i}-z_{j}\notin \hbar\CC[[\hbar]]$ for
$i\neq j$. 

\begin{prop} \label{conj}
Let $D^{(i)}$ be the element of $DY(\SL_{2})^{\otimes n}$ equal to 
$1^{\otimes (i-1)} \otimes D \otimes 1^{\otimes (n-i)}$. We have 
$$
\Ad(q^{\la D^{(i)}})(i_{\zz}( U_{\hbar} \G_{\zz})) 
= i_{\zz + \hbar
\la\de_{i}} ( U_{\hbar} \G_{\zz + \la \hbar \de_{i}}),
$$
for $\la\in \CC$ and $\de_{i}$ the $i$th standard basis vector of
$\CC^{n}$. 
\end{prop}

{\em Proof.}
$i_{\zz}( U_{\hbar} \G_{\zz})$ is generated by the 
coefficients of (\ref{isom:L-}). We have for any $j = 1,\ldots,n$, 
$$
\Ad(q^{\la D^{(i)}})(i_{\zz}(L_{(j)}^{-}(t_{j})))
= i_{\zz+\hbar \la\de_{i}} (L_{(j)}^{-}(t_{j}-\hbar \la\de_{ij})); 
$$
these are generating functions for $i_{\zz +
\hbar\la\de_{i}}(U_{\hbar}\G_{\zz+\hbar\la\de_{i}})$. 
\hfill \qed

\section{Quantum deformation of the zero-mode of the Sugawara field}

\subsection{Quantum Casimir elements}

In \cite{Dr}, V. Drinfeld proved the following fact: 

\begin{prop} (see \cite{Dr}, Prop. 2.2) \label{dr}
Let $\cA$ be a quasitriangular 
Hopf algebra with coproduct $\Delta_{\cA}$, counit
$\varepsilon_{\cA}$ and antipode $S_{\cA}$. Let $\cR_{\cA}$ be its
$R$-matrix and set 
$\cR_{\cA}^{-1} 
= \sum_{i} c_{i} \otimes d_{i}$. Let
$$
u = \sum_{i} d_{i} S_{\cA}(c_{i}).  
$$
Then we have for any $x$ in $\cA$, $S_{\cA}^{-2}(x) = uxu^{-1}$. 
\end{prop}

In particular, if for some other $u_{0}\in \cA$, we have $S_{\cA}^{-2}(x)
= u_{0} x u_{0}^{-1}$, then $u_{0}^{-1}u$ belongs to the center of
$\cA$. 

\subsection{Application to the deformation of the zero-mode of the
Sugawara tensor}

\label{exYg}

Let $\bar\G = \SL_{2}$ and 
$$
\G = \left(\bar\G \otimes \CC((t)) \right) \oplus \CC D \oplus \CC K 
$$
be the double extension of the loop algebra $\bar\G \otimes
\CC((t))$ by the cocycle 
$$
c(x\otimes \phi , y\otimes \psi) = \langle x,y\rangle_{\bar\G} \res_{0}(d\phi
\psi) K, 
$$
($\langle , \rangle_{\bar\G}$ is the Killing form of $\bar\G$), 
and by the derivation $[D,x\otimes \phi] = x \otimes ({d\phi}/{dt})$. 

Let $\wt{U \G}$ be the completion of the enveloping algebra
of $\G$, defined by the left ideals generated by the $\bar\G \otimes
t^{N}\CC[[t]]$. 

Let $e,f,h$ be the Chevalley basis of $\SL_{2}$, and let us set $x_{n} =
x\otimes t^{n}$ for $x \in \SL_{2}$. Then it is a known fact that
$$
(K+2)D + \sum_{k\ge 0}e_{-k-1}f_{k} + f_{-k-1}e_{k} + {1\over
2}h_{-k-1}h_{k} 
$$
belongs to the center of $\wt{U\G}$. The sum in this expression is the
zero-mode of the Sugawara tensor. 

Let $DY(\SL_{2})$ be the double Yangian algebra associated with
$\SL_{2}$. As we have seen, 
this algebra is generated by central and derivation elements
$K$ and $D$, and elements $\wt x$ lifting elements $x$ of $\bar\G
\otimes \CC((t))$.
Let $DY(\SL_{2})'$ be the subalgebra of $DY(\SL_{2})$
generated by $K$ and the $\wt x$. The universal $R$-matrix of
$DY(\SL_{2})$ is expressed as
$$
\cR_{\Yg} = q^{D\otimes K} \cR^{0}_{\Yg}, \quad \cR^{0}_{\Yg} \in
(DY(\SL_{2})')^{\bar\otimes 2}
$$
(we set $q = e^{\hbar}$). 
Moreover, $\cR^{0}_{\Yg}$ has the expansion 
$$
\cR^{0}_{\Yg} 
= 1 + \hbar \sum_{k\ge 0}\wt e_{k} \otimes \wt f_{-k-1} + \wt
f_{k} \otimes \wt e_{-k-1} + {1\over
2}\wt h_{k} \otimes \wt h_{-k-1} + O(\hbar^{2}). 
$$
The antipode $S_{\Yg}$ of $DY(\SL_{2})$ satisfies $S_{\Yg}^{2} = 
\Ad(q^{2D})$. 

Let us set 
$$
(\cR_{\Yg}^{0})^{-1} = \sum_{i} c_{i}^{0} \otimes d_{i}^{0}, 
$$
$c_{i}^{0}, d_{i}^{0}\in DY(\SL_{2})'$. 

\begin{prop} \label{Sug}
The sum 
$$
T = \sum_{i} q^{2D} d_{i}^{0} q^{KD} S_{\Yg}(c_{i}^{0})
$$
is a central element of
$DY(\SL_2)$. It is written as 
\begin{equation} \label{f:T}
T = q^{(K+2)D} S = S q^{(K+2)D}, 
\end{equation}
where $S$ belongs to the completion of the subalgebra 
$DY(\SL_{2})'$ defined by the left ideals 
generated by the lifts of $\bar\G \otimes t^{N}\CC[[t]]$. 
Its expansion in powers of $\hbar$ is 
\begin{equation} \label{expan:T}
T = 1 + \hbar \left((K+2)D +
\sum_{k\ge 0}\wt e_{-k-1}\wt f_{k} + \wt f_{-k-1}\wt e_{k} + {1\over
2}\wt h_{-k-1}\wt h_{k}  \right) + O(\hbar^{2}). 
\end{equation}
\end{prop}

{\em Proof.}
The first statement follows from Prop. \ref{dr}. We then have
$S = \sum_{i} (q^{-KD} d_{i}^{0} q^{KD})S_{\Yg}(c_{i}^{0})$. Since 
$\cR^{0}_{\Yg}$ commutes with $D\otimes 1 + 1\otimes D$, 
$S$ commutes with $D$. 
This proves (\ref{f:T}). 
(\ref{expan:T}) then follows from the above expansion of $\cR_{0}$. 
\hfill\qed\medskip

\section{Discrete connection on coinvariants}

\subsection{Induced representations}

Let $DY(\SL_{2})_{\ge 0}$ be the subalgebra of $DY(\SL_{2})$
generated by the the nonnegative Fourier generators $x_k,k\ge 0$, $x =
e,f,h$.  $DY(\SL_{2})_{\ge 0}^{\otimes n}$ 
is isomorphic to its image 
in $DY(\SL_{2})^{\otimes n} / (K^{(i)} - K^{(j)})$; these two algebras
will be denoted the same way. Finally, we denote by 
$DY(\SL_{2})_{\ge 0}^{\otimes n}[K]$ the subalgebra of
$DY(\SL_{2})^{\otimes n} / (K^{(i)} - K^{(j)})$ generated by 
$DY(\SL_{2})_{\ge 0}^{\otimes n}$ and $K$. 
We then have: 

\begin{lemma}
$i_{\zz}$ restricts to an isomorphism between $U_{\hbar}\G_{\cO}$ 
and $DY(\SL_{2})_{\ge 0}^{\otimes n}[K]$. 
\end{lemma} 

{\em Proof.}
It follows from (\ref{isom:L+}) and (\ref{rest}) that 
$i_{\zz}(U_{\hbar}\G_{\cO})$ is contained in 
$DY(\SL_{2})_{\ge 0}^{\otimes n}[K]$. Since the classical
limits of both algebras coincide with $U\G_{\cO}$ and the classical 
limit of $i_{\zz}$ is then the identity, $i_{\zz}$ is an isomorphism 
between these algebras.  
\hfill \qed \medskip

The restriction to $DY(\SL_{2})_{\ge 0}$  of the Yangian version of
$\pi$ can be specialized to $t=0$. Denote by $(V,\rho_{V})$ the
resulting $2$-dimensional representation. 

\begin{lemma}
We have
$$
( id \otimes \rho_{V}
) (L^{+}(t)) = R(t) .  
$$
\end{lemma}

{\em Proof.} We have 
$$
(id \otimes \rho_{V})(L^{+}(t)) = (\pi_{t} \otimes
\rho_{V})(\cR^{(21)}_{\Yg}) 
= (R^{-}(-t))^{(21)} = R(t), 
$$
where the third equality follows from (\ref{im:R}). \hfill \qed\medskip

We will consider also the dual representation 
                                                 $\rho_{V^{*}}$ defined on
$V^{*}$ by 
 $\rho_{V^{*}}(x) = \rho_{V}(S_{\Yg}^{-1}(x))^{t}$, where $t$
denotes the transposition. 
  
We have then 

\begin{lemma}
$$
( id \otimes \rho_{V^{*}} ) 
(L^{+}(t)) = (R(t)^{t_{2}})^{-1},  
$$
where the exponent $t_{2}$ denotes the transposition in the
second factor.
\end{lemma}  

{\em Proof.} We have (see \cite{Dr}, Prop. 3.1)
$(S_{\Yg} \otimes 1)(\cR_{\Yg}) = \cR_{\Yg}^{-1}$, so that $(1\otimes
\rho_{V})(\Id\otimes
S_{\Yg})(L^{+}(t)) = R(t)^{-1}$. On the other hand, $(\Id \otimes
S_{\Yg}^{-2})(L^{+}(t)) = L^{+}(t+2\hbar)$. Therefore,  
$(\Id \otimes \rho_{V})(\Id\otimes S_{\Yg}^{-1})(L^{+}(t)) =
R(t+2\hbar)^{-1}.$
 So we have 
$$
(\Id \otimes \rho_{V}^{t})(\Id\otimes S_{\Yg}^{-1})(L^{+}(t)) =
(R(t+2\hbar)^{-1})^{t_{2}}.
$$ 
The result now follows from the identity 
$$
(R(t+2\hbar)^{-1})^{t_{2}}
= 
(R(t)^{t_{2}})^{-1}.
$$
\hfill \qed\medskip

Let us fix now a complex number $k$ and consider the module 
$(\rho,(V^*)^{\otimes n})$ 
over $DY(\SL_{2})_{\ge 0}^{\otimes n}[K]$, defined 
as follows. As an algebra, $DY(\SL_{2})_{\ge 0}^{\otimes n}[K]$ is 
isomorphic to $DY(\SL_{2})_{\ge 0}^{\otimes n} 
\otimes \CC[K]$.  $DY(\SL_{2})_{\ge 0}^{\otimes n}$ then acts
on $(V^*)^{\otimes n}$ by $\rho_{V^{*}}^{\otimes n}$, and $K$ acts on this space 
by the scalar $k$. 

Let $\VV$ be the induced module 
$$
\VV = ( DY(\SL_2)' )^{\otimes n} / (K^{(i)} - K^{(j)})
\otimes_{DY(\SL_{2})_{\ge 0}^{\otimes n}[K]}
(V^*)^{\otimes n} ; 
$$
this is a module over the algebra 
$( DY(\SL_2)')^{\otimes n} / (K^{(i)} - K^{(j)})$. 

\subsection{Coinvariants}

We will be interested in 
the (dual to the) space of coinvariants
$$
H_{0}(U_{\hbar}\G_{\zz}, \VV)^{*},
$$ 
which is defined as the subspace of
$\VV^{*}$ consisting of the forms $\ell_{\zz}$, such that
$\ell_{\zz}(i_{\zz}(x)v) = \varepsilon(x)\ell_{\zz}(v)$, 
for $x\in U_{\hbar}\G_{\zz}$. 

Props. \ref{PBW} and \ref{isom} show that the map from 
$(V^*)^{\otimes n}$ to $\VV$, sending $v$ to $1\otimes v$, is injective. 
This map will serve to identify $(V^*)^{\otimes n}$ with a subspace of $\VV$. 

\begin{prop} \label{isos}
The restriction of $\ell_{\zz}$ to $(V^*)^{\otimes n}$ defines a map from
$H_{0}(U_{\hbar}\G_{\zz}, \VV)^{*}$ 
to $((V^*)^{\otimes n})^{*}=V^{\otimes n}$. This map
is a linear isomorphism. 
\end{prop}

{\em Proof.} It follows from Props. \ref{isom} and \ref{PBW} 
that a basis of $\VV$ is given by the $u_i \otimes v_j$, with 
$(u_i),(v_j)$ bases of $U_{\hbar}\G_{\zz}$ 
and $(V^*)^{\otimes n}$ respectively ($u_0 = 1$). A form $\ell_{\zz}$
of $\VV^*$ is then invariant iff it satisfies the equations
$\ell_{\zz} (u_i \otimes v_j) 
= \varepsilon(u_i)\ell_\zz(v_j)$, 
so it is exactly determined by the $\ell_{\zz} (1 \otimes v_j)$. 
\hfill \qed

\subsection{Compatible difference system on coinvariants} \label{cnx}

\subsubsection{Definitions}

Let $\CC^{n}_{*}$ be the complement of the diagonals in $\CC^{n}$, 
and let $\eta$ be a nonzero complex number. 
A difference flat connection on $\CC^{n}_{*}$ is the data of
a vector space $E_{\zz}$ for each $\zz\in \CC^{n}_{*}$, together with 
a system of linear isomorphisms 
$$
A_{i}(\zz): E_{\zz} \to E_{\zz + \eta\de_{i}}, 
$$
satisfying the relations
\begin{equation} \label{compat}
A_{j}(\zz + \eta\de_{i}) \circ A_{i}(\zz)
=
A_{i}(\zz + \eta\de_{j}) \circ A_{j}(\zz), \quad
\forall i,j = 1,\ldots,n; 
\end{equation}
$\eta$ is called the step of the system. 

Suppose we have a system $\iota_{\zz}$ of
isomorphisms of the $E_{\zz}$ with a fixed vector space $E$. 
We get a system of elements $\wt A_{i}(\zz) = \iota_{\zz +\hbar\de_i} 
A_i(\zz)\iota_{\zz}^ {-1}\in \Aut(E)$, satisfying 
the  same relations. Such a system is called a compatible
difference system (see \cite{A,TV}). 

\subsubsection{The case of a formal step}

We will need the following modification of the above definitions
in the formal context. In that situation, $\zz$ is a sequence 
$(z_{1},\ldots,z_{n}) \in \CC[[\hbar]]^{n}$, such that $z_{i}\neq z_{j}$ 
mod $\hbar$ for $i\neq j$.  The $E_{\zz}$ are free $\CC[[\hbar]]$-modules
and $\eta$ belongs to $\hbar\CC[[\hbar]]$. 

A compatible difference system will be a system of 
$\CC[[\hbar]]$-linear
isomorphisms $A_{i}(\zz): E_{\zz} \to E_{\zz + \eta\de_{i}}$ satisfying
(\ref{compat}); a system of $\CC[[\hbar]]$-linear 
isomorphisms $\iota_{\zz}:E_{\zz} \to E$ then gives
rise to a system of elements $\wt A_{i}\in \Aut_{\CC[[\hbar]]}(E)$, 
satisfying the same relations, again called a compatible difference 
system. (In the case where $E$ has rank $p$, $\Aut_{\CC[[\hbar]]}(E)$
is isomorphic to $GL_{p}(\CC[[\hbar]])$.)

\subsubsection{The action of the quantum Sugawara element on
coinvariants} \label{act:Sug}

We then have a compatible difference system, in the formal sense, 
defined as follows. Associate to $\zz$ the coinvariants space
$E_{\zz} = H_{0}(U_{\hbar}\G_{\zz} , \VV)^{*}$. Set 
$\eta = \hbar(k+2)$ and define 
$$
A_{i}(\zz): E_{\zz} \to E_{\zz + \hbar(k+2)\de_{i}}
$$
by the formula 
\begin{equation} \label{conn}
(A_{i}(\zz)\ell_{\zz})(v) = \ell_{\zz}(\rho(S^{(i)})v), 
\end{equation}
for $\ell_{\zz}\in E_{\zz}, v\in \VV_{\la}$, where 
$1^{\otimes (i-1)} \otimes S \otimes 1^{\otimes (n-i)}$ 
and $S$ is the quantum Sugawara element 
defined in Prop. \ref{Sug}. 

Let us show that the form $A_{i}(\zz)\ell_{\zz}$ in invariant with respect 
to 
$
U_{\hbar}\G'_{\cK,\zz + \hbar(k+2)\de_{i}}.
$
For $x\in U_{\hbar}\G'_{\cK,\zz + \hbar(k+2)\de_{i}}$, we have
\begin{align*}
(A_{i}(\zz)\ell_{\zz})( i_{\zz + \hbar(k+2)\de_{i}}(x) v) & = 
\ell_{\zz}(S^{(i)} i_{\zz + \hbar(k+2)\de_{i}}(x)v)
\\ & 
= \ell_{\zz}( \Ad(S^{(i)})( i_{\zz + \hbar(k+2)\de_{i}}(x)) S^{(i)}v) ; 
\end{align*}
but since $q^{(K+2)D}S$ is central in $DY(\SL_{2})$, 
we have $\Ad(S^{(i)}) (i_{\zz + \hbar(k+2)\de_{i}}(x))
= \Ad(q^{ - (K+2) D^{(i)}}) (i_{\zz + \hbar(k+2)\de_{i}}(x))$; 
the action of this element on $\VV$ coincides with that of 
$\Ad(q^{ - (k+2) D^{(i)}}) (i_{\zz + \hbar(k+2)\de_{i}}(x))$, 
which belongs to  $i_{\zz} (U_{\hbar}\G_{\zz})$ by Prop.~\ref{conj}. 

We then have 
\begin{align*}
(A_{i}(\zz)\ell_{\zz})( i_{\zz + \hbar(k+2)\de_{i}}(x) v) & =
\ell_{\zz} ( \Ad(q^{ - (k+2) D^{(i)}}) (i_{\zz + \hbar(k+2)\de_{i}}(x)) 
v) \\ & 
= \varepsilon( \Ad(q^{ - (k+2) D^{(i)}}) (i_{\zz + \hbar(k+2)\de_{i}}(x)) ) 
\ell_{\zz}(v) \\ & 
= \varepsilon(x)\ell_{\zz}(v). 
\end{align*}

This shows that $A_{i}(\zz)\ell_{\zz}$ is $U_{\hbar}\G_{\zz 
+ \hbar(k+2)\de_{i}}$-invariant.  

\section{Identification with the qKZ system}

The aim of this section is to make the system (\ref{conn}) explicit, 
using the identifications  of Prop. \ref{isos} of the spaces of 
coinvariants with  $V^{\otimes n}$. 

\subsection{Expression of the Sugawara action in terms of $L$-operators}

Let us express the connection (\ref{conn}) in terms of $L$-operators. 
For that, we will prove: 

\begin{lemma} \label{exp:Sug}
Let us write $L^{-}(\hbar k) = \sum_{\al} a_{\al} \otimes l^{-}_{\al}$, 
with $a_{\al}\in \End(V)[[\hbar]]$ and $l^{-}_{\al}\in DY(\SL_{2})$.
Let $v$ belong to $(V^*)^{\otimes n}$; we view it as a vector of $\VV$, as 
explained above. Then the action of $S^{(i)}$ on $v$ is expressed as
$$
\rho(S^{(i)}) v = \sum_{\al} \rho((l^{-}_{\al})^{(i)}) 
(a_{\al}^{t})^{(i)} v ; 
$$
recall that
the exponent $i$
means the action on the $i$th factor of $(V^*)^{\otimes n}$. 
\end{lemma}

{\em Proof.}
In the notation of sect.\ \ref{exYg}, we have $S = \sum_{i} 
(q^{-KD} d_{i}^{0} q^{KD}) S_{\Yg}(c_{i}^{0})$. Therefore, 
\begin{equation} \label{ronco}
\rho(S^{(i)})  = \sum_{j} \rho(q^{-kD} d_{j}^{0} q^{kD})
(1^{\otimes (i-1)} \otimes \rho_{V^{*}}(S_{\Yg}(c_{j}^{0})) 
\otimes 1^{\otimes (n-i)}).  
\end{equation}

On the other hand, we have $(1\otimes q^{-kD}) 
L^{-}(t) (1\otimes q^{kD}) = L^{-}(t+\hbar k)$, 
so that for $t=0$ we obtain, since $\rho_{V}$
coincides with the specialization of $\pi$ for $t=0$,  
$$
\sum_{i} \rho_{V}(c_{j}^{0}) \otimes q^{-KD} d_{j}^{0} q^{KD} = 
\sum_{\al} a_{\al} \otimes l_{\al}^{-}. 
$$

Therefore, we have
$$
\sum_{i} \rho_{V^{*}}(S_{\Yg}(c_{j}^{0})) \otimes q^{-KD} d_{j}^{0} q^{KD} = 
\sum_{\al} a_{\al}^{t} \otimes l_{\al}^{-}. 
$$ 

The lemma follows from the comparison of this formula with 
(\ref{ronco}). 
\hfill \qed\medskip 

\subsection{Expression of the discrete compatible system}

Let us express the invariance of the form $\ell_{\zz}$. Let $V^{a}$ 
be an auxiliary copy of the vector space $V$. Consider 
$i_{\zz}(L^{-}_{(i)}(t_{i}))$ as an elements of $\End(V^{a}) \otimes 
DY(\SL_{2})^{\otimes n} / (K^{(i)} - K^{(j)})$. 
We then have, for $v_{a}\in 
V^{a}, v\in (V^*)^{\otimes n}$, and any formal $t_{i}$, 
$$
(1\otimes \ell_{\zz})\left((1\otimes \rho)(i_{\zz}(L^{-}_{(i)}(t_{i})))
(v_{a} \otimes v) \right) 
=
(1\otimes\ell_{\zz})(v_{a} \otimes v). 
$$
Substitute 
$t_{i}$ by $\hbar k$ in this identity. Using (\ref{isom:L-}), 
we find that 
\begin{align*}
\sum_{\al} (1\otimes \ell_{\zz}) & \left(  
( R^{(a1)}(\hbar k + z_{i,1})^{t_{1}} )^{-1} \cdots
( R^{(a,i-1)}(\hbar k + z_{i,i-1})^{t_{i-1}} )^{-1}
\rho(l_{\al}^{-(i)}) a_{\al}^{(a)}
\right. \\ & \left. 
( R^{(a,i+1)}(z_{i,i+1})^{t_{i+1}})^{-1} \cdots
( R^{(an)}(z_{i,n})^{t_{n}} )^{-1} 
(v^{a} \otimes v)\right) \\ & 
= (1\otimes \ell_{\zz})(v^{a} \otimes v). 
\end{align*}
The exponent $t_i$ denotes here the transposition of the
$i$th factor.

Therefore, we have the identity
\begin{align} \label{id}
& \sum_{\al} (1\otimes \ell_{\zz})  \left(  
( R^{(a1)}(\hbar k + z_{i,1})^{t_{1}})^{-1} \cdots
( R^{(a,i-1)}(\hbar k + z_{i,i-1})^{t_{i-1}})^{-1}
\rho(l_{\al}^{-(i)}) a_{\al}^{(a)}
\right. \\ & \left. \nonumber  
(v^{a} \otimes v)\right) \\ & \nonumber
= (1\otimes \ell_{\zz})
\left( 
R^{(a,n)}(z_{i,n})^{t_{n}} \cdots
R^{(a, i+1)}(z_{i,i+1})^{t_{i+1}} 
(v^{a} \otimes v)
\right) . 
\end{align}
for any $v^{a}\in V^{a}$, $v\in (V^*)^{\otimes n}$. 

Introduce now
\begin{equation*}
\tilde R(z)= \left( \left(R(z)^{-1}\right)^{t_2}\right)^{-1}
\end{equation*}
This object  has the following property:
if we write $(R(z)^{t_{2}})^{-1}=\sum_i x'_i\otimes x''_i$ and
$\tilde R(z)
                     =\sum_iy'_i\otimes y''_i$, then
$\sum_{ij}y'_ix'_j\otimes x''_j y''_i=\Id_{V\otimes V^*}$. Thus we
have
\begin{align}\label{id-2}
 & \sum_{\al} (1\otimes \ell_{\zz})
\left( \rho(l_{\al}^{-(i)}) a_{\al}^{(a)}
(v^{a} \otimes v) \right)  \\ & \nonumber
= (1\otimes \ell_{\zz})
\bigl(  
\tilde R^{(a,i-1)}(\hbar k + z_{i,i-1})\cdots
\tilde R^{(a1)}(\hbar k + z_{i,1}) 
\\
 &\
 R^{(a,n)}(z_{i,n})^{t_{n}} \cdots
R^{(a,i+1)}(z_{i,i+1})^{t_{i+1}} 
(v^{a} \otimes v)
\bigr). \nonumber
\end{align}

Note that we have $\tilde R(z) = 
                                      R(z+2\hbar)^{t_{2}}      $, so that 
\begin{align*}
 & \sum_{\al} (1\otimes \ell_{\zz})
\left( \rho(l_{\al}^{-(i)}) a_{\al}^{(a)}
(v^{a} \otimes v) \right)  \\ & \nonumber
= (1\otimes \ell_{\zz})
\left( R^{(a,i-1)}(z_{i,i-1} + \hbar(k+2))^{t_{i-1}} 
\cdots R^{(a,1)}(z_{i,1} + \hbar(k+2))^{t_{1}}  \right. \\ 
 &\
 R^{(a,n)}(z_{i,n})^{t_{n}} \cdots
R^{(a,i+1)}(z_{i,i+1})^{t_{i+1}} 
(v^{a} \otimes v)
\nonumber \left. \right) . 
\end{align*}

After we apply the transposition on $V^{a}$, we obtain,
for any $u^a\in V^*$,
\begin{align} \label{id-3}
 & \sum_{\al} (1\otimes \ell_{\zz})
\left( \rho(l_{\al}^{-(i)}) (a_{\al}^{(a)})^{t}
(u^{a} \otimes v) \right)  \\ & \nonumber
= (1\otimes \ell_{\zz})
\left( R^{(a,i+1)}(z_{i,i+1})^t
\cdots R^{(a,n)}(z_{i,n}) ^t
\right. \\ 
 &\  R^{(a,1)}(z_{i,1} + \hbar(k+2))^t  \cdots
R^{(a,i-1)}(z_{i,i-1} + \hbar(k+2))^t 
(u^{a} \otimes v)
\nonumber \left. \right) . 
\end{align}
Here $R(z)^t=R(z)^{t_1t_2}$ is the transposed of $R(z)$.

Let $v_{1},\ldots, v_{n}$ belong to $V^*$. Let $(e_{\beta})$ be a basis of 
$V^*$ 
and $(e^{\beta})$ be the dual basis. Let us apply the identity (\ref{id-3})
to 
$$
u^{a} \otimes v
= v_{i} \otimes (v_{1}\otimes\cdots \otimes v_{i-1} 
\otimes e_{\beta} \otimes v_{i+1} \otimes\cdots\otimes v_{n}), 
$$
apply $e^{\beta} \otimes 1$ to this identity, and sum over $\beta$. We find, 
with $v_{1,i-1} = v_{1} \otimes\cdots \otimes v_{i-1}$ 
and $v_{i+1,n} = v_{i+1} \otimes\cdots \otimes v_{n}$,  
\begin{align} \label{id-4}
& \sum_{\al}\sum_{\beta} (e^{\beta}\otimes \ell_{\zz})  
\left(  
\rho(l_{\al}^{-(i)}) (a_{\al}^{(a)})^{t} 
(v_{i} \otimes v_{1,i-1} \otimes e_{\beta} \otimes v_{i+1,n})
\right)
\\ & \nonumber
= \sum_{\beta} (e^{\beta}\otimes \ell_{\zz})
\left(  R^{(a,i+1)}(z_{i,i+1})
\cdots R^{(a,n)}(z_{i,n}) 
\right. 
\\ &\  R^{(a,1)}(z_{i,1} + \hbar(k+2))  \cdots
R^{(a,i-1)}(z_{i,i-1} + \hbar(k+2))
(v_{i} \otimes v_{1,i-1} \otimes e_{\beta} \otimes v_{i+1,n})
\left. \right) . 
\end{align}

The following result is a consequence of the definition of 
dual bases. 

\begin{lemma}
For any endomorphism $X$ of $V^*$, and any element $v$ of $V^*$, 
we have 
$$
\sum_{\al}(e^{\beta} \otimes 1)(Xv \otimes e_{\beta}) = Xv.
$$
\end{lemma}
Applying this lemma to (\ref{id-3}), we find 
\begin{align} \label{id4}
& \sum_{\al} \ell_{\zz}  (  
\rho(l_{\al}^{-(i)}) (a_{\al}^{(i)})^{t} v )
= \ell_{\zz}
\left( 
R^{(i,i+1)}(z_{i,i+1})^t \cdots
R^{(in)}(z_{i,n})^t 
\right.  \\ & \nonumber \left. 
R^{(i1)}(z_{i,1} + \hbar(k+2))^t \cdots
R^{(i,i-1)}(z_{i,i-1} + \hbar(k+2))^t v
\right) . 
\end{align}

By Lemma \ref{exp:Sug}, the left side of this equality is the expression of
$(A_{i}\ell_{\zz})(v)$. 

\begin{thm}
The discrete flat connection on the spaces of coinvariants defined by  
(\ref{conn}) is identified by the isomorphisms (\ref{isos}) with the 
quantum Knizhnik-Zamolodchikov system 
\begin{align}\label{QKZ}
A_{i}\ell_{\zz} = & 
R^{(i,i-1)}(z_{i,i-1} + \hbar(k+2))
\cdots R^{(i1)}(z_{i,1} + \hbar(k+2))  
 \\  
& \ \ \ \ \  
R^{(in)}(z_{i,n}) \cdots R^{(i,i+1)}(z_{i,i+1}) 
\ell_{\zz}.\nonumber 
\end{align}
with step $\eta = \hbar(k+2)$.
\end{thm}

\begin{remark} We have obtained equations for the coinvariants of the
representation $\VV$ of $U_{\hbar}\G_{\cK,\zz}$ induced from
$\rho_{V^{*}}^{\otimes n}$. It should be possible to obtain equations
for coinvariants of a representation $\otimes_{i=1}^{n}
\rho_{V_{i}^{*}}$, where $V_{i}$ are finite dimensional representations
of $DY(\SL_{2})_{\ge 0}$. For that, one should consider the $L$-operators
$L^{-(V_{i})}_{(i)}(t_{i})$ for $U_{\hbar}\G_{\cK,\zz}$ and 
$L^{\pm(V_{i})}(t)$, defined by replacing $\pi$ by $\pi_{V_{i}}$ in the
definitions of $L^{-}_{(i)}(t_{i})$ and $L^{\pm}(t)$, and prove that
$i_{\zz}(L^{-(V_{i})}_{(i)}(t_{i}))$ is given by the formula of
Prop. \ref{isos}. After that, it should be possible to apply the
reasoning above to derive the general version of the qKZ equation. 
\end{remark}

\begin{remark}
It is interesting to consider the systems obtained by replacing
$\VV$ by a representation induced from a non-irreducible
representation of $U_{\hbar}\G_{\cO}$. In the classical case and the
$\GL_{1}$ situation, one obtains in this way non-Fuchsian systems. 

For example, $DY(\SL_{2})_{\ge 0}$ has
an evaluation representation $\rho_{V[[\zeta]]}$ on $V[[\zeta]]$,
defined by $(id \otimes \rho_{V[[\zeta]]})(L^{+}(t)) = R(t+\zeta)$;  
the representation $\rho_{V}$ considered above is a 
quotient of $\rho_{V[[\zeta]]}$. 
One may expect that the system of equations one would obtain this way
is the system 
$$
\eta = \hbar(k+2), \quad
A_{i}\ell_{\zz} = K_{i}(z_{ij}+\zeta_{i}-\zeta_{j})
\ell_{\zz}, 
$$
where $\ell_{\zz}$ belongs to $V^{\otimes
n}[[\zeta_{1},\ldots,\zeta_{n}]]$, on which the operators
$K_{i}(z_{ij}+\zeta_{i}-\zeta_{j})$ act 
as $\sum_{k\ge 0}K_{i}^{(k)}(z_{ij})(\zeta_{i}-\zeta_{j})^{k}/k!$
(the functions of $\zeta_{i}$ act by multiplication on the formal
series part and the
derivatives of $K_{i}$ act as matrices on the factor $V^{\otimes
n}$; $K_{i}$ are the operators appearing in the right side of
(\ref{QKZ})).  
\end{remark}

\begin{remark} In the classical case, it is possible to explain the
agreement of the intertwiners and coinvariants approaches in a simple
way. It would be interesting to find such an explanation of the result
of this paper. 
\end{remark}

\end{document}